\input amstex
\documentstyle{amsppt}
\magnification=\magstep1
 \hsize 13cm \vsize 18.35cm \pageno=1
\loadbold \loadmsam
    \loadmsbm
    \UseAMSsymbols
\topmatter
\NoRunningHeads
\title Symmetry $p$-adic  invariant integral on $\Bbb Z_p$ for
Bernoulli and Euler polynomials
\endtitle
\author
  Taekyun Kim 
\endauthor
 \keywords $q$-Bernoulli numbers, $q$-Volkenborn integrals,
 $q$-Euler numbers, $q$-Stirling numbers
\endkeywords

\abstract The main purpose of this paper is to investigate several
further interesting properties of symmetry for the $p$-adic
invariant integrals on $\Bbb Z_p$. From these symmetry , we can
derive many interesting recurrence identities for Bernoulli and
Euler polynomials. Finally we introduce the new concept of symmetry
of fermionic $p$-adic invariant integral on $\Bbb Z_p$. By using
this symmetry of fermionic $p$-adic invariant integral on $\Bbb
Z_p$, we will give some relations of symmetry between the power sum
polynomials and Euler numbers. The relation between the
$q$-Bernoulli polynomials and $q$-Dedekind type sums which discussed
in [ Y. Simsek, q-Dedekind type sums related to q-zeta function and
basic L-series, J. Math. Anal. Appl.318(2006), 333-351] can be also
derived by using the properties of symmetry of fermionic $p$-adic
integral on $\Bbb Z_p$.

\endabstract
\thanks  2000 AMS Subject Classification: 11B68, 11S80
\newline  This paper is supported by  Jangjeon Research Institute for Mathematical
Science(JRIMS-10R-2001)
\endthanks
\endtopmatter

\document

{\bf\centerline {\S 1. Introduction}}

 \vskip 20pt

Let $p$ be a fixed prime. Throughout this paper $\Bbb Z_p ,$ $\Bbb
Q_p ,$ $\Bbb C,$ and $\Bbb C_p$ will, respectively, denote the ring
of $p$-adic rational integers, the field of $p$-adic rational
numbers, the complex number field, and the completion of algebraic
closure of $\Bbb Q_p .$  For $x\in\Bbb C_p ,$ we use the notation
$[x]_q=\frac{1-q^x}{1-q}, $ cf. [1-6]. The Carlitz's $q$-Bernoulli
numbers $\beta_{k,q}$ can be determined inductively by
$$\eqalignno{ &
\beta_{0,q} =1,\quad  q(q\beta +1)^k -\beta_{k,q} = \cases
1 & \text{if\   $k=1$}\\
0 & \text{if \  $k>1$,}
\endcases }
$$
with the usual convention of replacing $\beta^i$ by
$\beta_{i,q}$(see [2, 3, 24, 25]).

We say that $f$ is a uniformly differentiable function at a point $a
\in\Bbb Z_p $ and denote this property by $f\in UD(\Bbb Z_p )$, if
the difference quotients $F_f (x,y) = \dfrac{f(x) -f(y)}{x-y} $ have
a limit $l=f^\prime (a)$ as $(x,y) \to (a,a)$. For $f\in UD(\Bbb Z_p
)$, let us start with the expression

$$\eqalignno{ & \dfrac{1}{[p^N ]_q} \sum_{0\leq j < p^N} q^j f(j) =\sum_{0\leq j < p^N} f(j)
\mu_q (j +p^N \Bbb Z_p ), }
$$
representing a $q$-analogue of Riemann sums for $f$, cf. [8-14]. The
integral of $f$ on $\Bbb Z_p$ will be defined as limit ($n \to
\infty$) of those sums, when it exists. The $q$-deformed bosonic
$p$-adic integral of the function $f\in UD(\Bbb Z_p )$ is defined by
$$ I_q (f)= \int_{\Bbb Z_p }f(x) d\mu_q (x) = \lim_{N\to \infty}
\dfrac{1}{[dp^N ]_q} \sum_{0\leq x < dp^N} f(x) q^x, \text{ see [8]}
. \tag1$$ From(1) we note that
$$I_{1}(f)= \lim_{q\rightarrow 1}I_q(f)=\int_{\Bbb
Z_p}f(x)dx=\lim_{N\rightarrow
\infty}\frac{1}{p^N}\sum_{x=0}^{p^N-1}f(x). \tag2$$ In [8], it was
shown that the Carlitz's $q$-Bernoulli numbers can be represented by
$p$-adic $q$-integral on $\Bbb Z_p$ as follows:
$$\int_{\Bbb
Z_p}[x]_q^md\mu_q(x)=\beta_{m,q}, \text{ $m\in\Bbb Z_{+}$.}$$ Thus
we have $\lim_{q\rightarrow 1}\int_{\Bbb Z_p}[x]_q^n
d\mu_q(x)=\int_{\Bbb Z_p}x^n dx=B_n$, where $B_n$ are the $n$-th
ordinary Bernoulli numbers. In the recent the $p$-adic invariant
integral on $\Bbb Z_p$ are studied by several researchers in the
area of  number theory and mathematical physics (see [1-25]). In
[23], by using q-Volkenborn integral(=$p$-adic invariant integral on
$\Bbb Z_p$), Y. Simsek constructed new generating functions of the
new twisted (h,q)-Bernoulli polynomials and numbers. By applying the
Mellin transformation to these generating functions, he also obtain
integral representations of the new twisted (h,q)-zeta function and
twisted (h,q)-L-function, which interpolate the twisted
(h,q)-Bernoulli numbers and generalized twisted (h,q)-Bernoulli
numbers at non-positive integers, respectively. The main purpose of
this paper is to investigate several further interesting properties
of symmetry for the $p$-adic invariant integrals on $\Bbb Z_p$. From
these symmetry , we can derive many interesting recurrence
identities for Bernoulli and Euler polynomials. Finally we introduce
the new concept of symmetry of fermionic $p$-adic invariant integral
on $\Bbb Z_p$. By using this symmetry of fermionic $p$-adic
invariant integral on $\Bbb Z_p$, we will give some relations of
symmetry between the power sum polynomials and Euler numbers. The
relation between the $q$-Bernoulli polynomials and $q$-Dedekind type
sums which discussed in [25] can be also derived by using the
properties of symmetry of fermionic $p$-adic integral on $\Bbb Z_p$

\vskip 20pt

{\bf\centerline {\S 2. Symmetry  $p$-adic invariant integral on
$\Bbb Z_p$}} \vskip 10pt

In this section we assume that $q\in\Bbb C_p$ with $|1-q|_p<1$.
 By (2), we
easily see that
$$\int_{\Bbb Z_p}f(x+1)dx=\int_{\Bbb Z_p}f(x)dx+f^{\prime}(0), \text{ where $f^{\prime}(0)=\frac{df(x)}{dx}|_{x=0}$}.\tag3 $$
Let $a$ be an integer. From (3) we derive
$$\int_{\Bbb Z_p}e^{atx}dx=\sum_{n=0}^{\infty}B_n\frac{a^nt^n}{n!},
\tag 4$$ where $B_n$ are the $n$-th ordinary Bernoulli numbers. Let
$B_n(x)$ be the $n$-th Bernoulli polynomials. Then we have
$$\sum_{n=0}^{\infty}\frac{B_n(x)}{n!}a^nt^n=\int_{\Bbb
Z_p}e^{(ay+x)t}dy=\sum_{n=0}^{\infty}\int_{\Bbb
Z_p}(ay+x)^ndy\frac{t^n}{n!}. \tag5$$ By (4) and (5), we see that
$$\sum_{n=0}^{\infty}B_n(x)\frac{a^nt^n}{n!}=\left(\sum_{n=0}^{\infty}\frac{x^n}{n!}t^n\right)
\left(\sum_{k=0}^{\infty}a^kB_k\frac{t^k}{k!}\right)=\sum_{n=0}^{\infty}\left(\sum_{k=0}^n
\binom{n}{k}B_ka^kx^{n-k}\right)\frac{t^n}{n!}.$$ Thus, we have
$$B_n(x)=\sum_{k=0}^n\binom{n}{k}B_k (\frac{x}{a})^{n-k} .$$
 From the iterative method, we derive
 $$\int_{\Bbb Z_p}f(x+n)dx=\int_{\Bbb
 Z_p}f(x)dx+\sum_{i=0}^{n-1}f^{\prime}(i), \text{ where $n\in\Bbb N$
 and  $f^{\prime}(i)=\frac{df(x)}{dx}|_{x=i}.$} \tag6$$
By (6), we have
$$\frac{1}{t}\left(\int_{\Bbb Z_p}e^{(n+x)t}dx-\int_{\Bbb
Z_p}e^{xt}dx \right)=\sum_{i=0}^{n-1}e^{it}=\sum_{k=0}^{\infty}
\left(\sum_{i=0}^{n-1}i^k \right)\frac{t^k}{k!}. \tag7$$ It is easy
to show that
$$\frac{1}{t}\left(\int_{\Bbb Z_p}e^{(n+x)t}dx-\int_{\Bbb
Z_p}e^{xt}dx\right)=\frac{n\int_{\Bbb Z_p}e^{xt}dx}{\int_{\Bbb
Z_p}e^{nxt}dx }.\tag8$$ For each integer $k\geq 0$, let $
S_k(n)=0^k+1^k+\cdots+n^k$. It is known that $S_k(n)$ is called the
sums of powers of consecutive integers. From (7) and (8), we note
that
$$\frac{1}{t}\left(\int_{\Bbb Z_p}e^{(n+x)t}dx-\int_{\Bbb
Z_p}e^{xt}dx \right)=\frac{n\int_{\Bbb Z_p}e^{xt}dx}{\int_{\Bbb
Z_p}e^{nxt}dx} =\sum_{k=0}^{\infty}S_k(n-1)\frac{t^k}{k!}. \tag9$$
In (9), we note that
$$\int_{\Bbb Z_p}(n+x)^k dx-\int_{\Bbb Z_p}x^ kdx=kS_{k-1}(n-1),
 \text{ $k \in\Bbb N $ .}\tag10 $$
The Eq.(10) is equivalent to the following Eq.(11).
$$\frac{ \left(B_k(n)-B_k
\right)}{k}=S_{k-1}(n-1)=\sum_{l=0}^{n-1}l^{k-1}, \text{where
$n,k\in\Bbb N .$}\tag11$$ Let $w_1, w_2 \in\Bbb N$. Then we easily
see that
$$\frac{\int_{\Bbb Z_p}\int_{\Bbb Z_p}e^{(w_1x_1+w_2x_2)t}dx_1
dx_2}{\int_{\Bbb
Z_p}e^{w_1w_2xt}dx}=\frac{t(e^{w_1w_2t}-1)}{(e^{w_1t}-1)(e^{w_2t}-1)}.\tag12$$
Now we also consider the following double $p$-adic invariant
integral on $\Bbb Z_p$ as follows:

$$I=\frac{\int_{\Bbb Z_p}\int_{\Bbb Z_p}e^{(w_1x_1+w_2x_2+w_1w_2x)t}dx_1
dx_2}{\int_{\Bbb
Z_p}e^{w_1w_2x_{3}t}dx_{3}}=\frac{te^{w_1w_2xt}(e^{w_1w_2t}-1)}{(e^{w_1t}-1)(e^{w_2t}-1)}.\tag13$$
From (7) and(8) we note that
$$\frac{w_1\int_{\Bbb Z_p}e^{xt}dx}{\int_{\Bbb
Z_p}e^{w_1xt}dx}=\sum_{l=0}^{\infty}\left(\sum_{k=0}^{w_1-1}k^l\right)\frac{t^l}{l!}
=\sum_{l=0}^{\infty}S_l(w_1-1)\frac{t^l}{l!}. \tag14$$ In (13) it is
easy to show that
$$I=\left(\frac{1}{w_1}\int_{\Bbb Z_p}e^{w_1(x_1+w_2x)t}dx_1 \right) \left(\frac{w_1\int_{\Bbb Z_p}e^{w_2x_2t}dx_2}
{\int_{\Bbb Z_p}e^{w_1w_2xt}dx}\right). \tag15$$ From (5), (14) and
(15), we can derive the following (16).
$$\aligned
I&=\frac{1}{w_1}\left(\sum_{i=0}^{\infty}B_i(w_2x)\frac{w_1^it^i}{i!}\right)\left(
\sum_{l=0}^{\infty}S_l(w_1-1)\frac{w_2^lt^l}{l!}\right)\\
&=\sum_{n=0}^{\infty}\left(\sum_{i=0}^n
B_i(w_2x)S_{n-i}(w_1-1)w_1^{i-1}w_2^{n-i}\right)t^n\\
&=\sum_{n=0}^{\infty}\left(\sum_{i=0}^n
\binom{n}{i}B_i(w_2x)S_{n-i}(w_1-1)w_1^{i-1}w_2^{n-i}\right)\frac{t^n}{n!}.
  \endaligned\tag16$$
By the symmetry of $p$-adic invariant integral on $\Bbb Z_p$, we
also see that
$$\aligned
I&=\left(\frac{1}{w_2}\int_{\Bbb
Z_p}e^{w_2(x_2+w_1x)t}dx_2\right)\left(\frac{w_2 \int_{\Bbb
Z_p}e^{w_1x_1t}dx_1}{\int_{\Bbb Z_p}e^{w_1w_2xt}dx} \right)\\
&=\frac{1}{w_2}\left(\sum_{i=0}^{\infty}B_i(w_1x)\frac{w_2^it^i}{i!}
\right)\left(\sum_{l=0}^{\infty}S_l(w_2-1)\frac{w_1^lt^l}{l!}\right)\\
&=\sum_{n=0}^{\infty}\left(\sum_{i=0}^n\binom{n}{i}B_i(w_1x)S_{n-i}(w_2-1)w_2^{i-1}w_1^{n-i}\right)\frac{t^n}{n!}.
\endaligned\tag17$$
By comparing the coefficients on the both sides of (16) and (17), we
obtain the following theorem and corollary:

\proclaim{Theorem 1} For all $w_1, w_2\in\Bbb N,$ we have
$$\aligned
&\left(\frac{1}{w_1}\int_{\Bbb Z_p}e^{w_1(x_1+w_2x)t}dx_1 \right)
\left(\frac{w_1\int_{\Bbb Z_p}e^{w_2x_2t}dx_2} {\int_{\Bbb
Z_p}e^{w_1w_2xt}dx}\right)\\
&=\left(\frac{1}{w_2}\int_{\Bbb
Z_p}e^{w_2(x_2+w_1x)t}dx_2\right)\left(\frac{w_2 \int_{\Bbb
Z_p}e^{w_1x_1t}dx_1}{\int_{\Bbb Z_p}e^{w_1w_2xt}dx}.
\right)\endaligned$$
\endproclaim

\proclaim{ Corollary 2} For $n\geq 0,$ we have
$$\sum_{i=0}^n
\binom{n}{i}B_i(w_2x)S_{n-i}(w_1-1)w_1^{i-1}w_2^{n-i}
=\sum_{i=0}^n\binom{n}{i}B_i(w_1x)S_{n-i}(w_2-1)w_2^{i-1}w_1^{n-i},$$
where $\binom{n}{i}$ is the binomial coefficients.
\endproclaim
If we take $w_2=1$ in Corollary 2, then we have
$$B_n(w_1x)=\sum_{i=0}^{n}w_1^{i-1}\binom{n}{i}B_i(x)S_{n-i}(w_1-1).$$
Thus, we can derive the formula of Deeba-Rodriguez as follows ( see
[26, 27]):
$$B_n=\frac{1}{w_1(1-w_1^n)}\sum_{k=0}^{n-1}w_1^k\binom{n}{k}B_kS_{n-k}(w_1-1).$$
By (5), (7) and (8), we see that
$$\aligned
I&=\left(\frac{e^{w_1w_2xt}}{w_1}\int_{\Bbb
Z_p}e^{w_1x_1t}dx_1\right)\left(\frac{w_1\int_{\Bbb
Z_p}e^{w_2x_2t}dx_2}{\int_{\Bbb Z_p}e^{w_1w_2xt}dx}\right)\\
&=\left(\frac{e^{w_1w_2xt}}{w_1}\int_{\Bbb
Z_p}e^{w_1x_1t}dx_1\right)\left(\sum_{i=0}^{w_1-1}e^{w_2it}\right)\\
&=\left(\frac{1}{w_1}\int_{\Bbb
Z_p}e^{w_1x_1t}dx_1\right)\left(\sum_{i=0}^{w_1-1}e^{(w_2x+\frac{w_2}{w_1}i)w_1t}\right)\\
&=\frac{1}{w_1}\sum_{i=0}^{w_1-1}\int_{\Bbb
Z_p}e^{(x_1+w_2x+\frac{w_2}{w_1}i)tw_1}dx_1
=\sum_{n=0}^{\infty}\left(\sum_{i=0}^{w_1-1}B_n(w_2x+\frac{w_2}{w_1}i)w_1^{n-1}\right)\frac{t^n}{n!}.
\endaligned\tag18$$
On the other hand, we obtain the following equation by the symmetry
of $p$-adic invariant integral on $\Bbb Z_p$ as follows:
$$\aligned
I&=\left(\frac{e^{w_1w_2xt}}{w_2}\int_{\Bbb
Z_p}e^{w_2x_2t}dx_2\right)\left(\frac{w_2\int_{\Bbb
Z_p}e^{w_1x_1t}dx_1}{\int_{\Bbb Z_p}e^{w_1w_2xt}dx}\right)\\
&=\left(\frac{1}{w_2}\int_{\Bbb
Z_p}e^{w_2x_2t}dx_2\right)\left(\sum_{i=0}^{w_2-1}e^{(w_1x+\frac{w_1}{w_2}i)w_2t}\right)\\
&=\frac{1}{w_2}\sum_{i=0}^{w_2-1}\int_{\Bbb
Z_p}e^{(x_2+w_1x+\frac{w_1}{w_2}i)tw_2}dx_2
=\sum_{n=0}^{\infty}\left(\sum_{i=0}^{w_2-1}B_n(w_1x+\frac{w_1}{w_2}i)w_2^{n-1}\right)\frac{t^n}{n!}.
\endaligned\tag18-1$$
Therefore we obtain the following symmetry for the $p$-adic
invariant integral on $\Bbb Z_p$ as follows:

\proclaim{Theorem 3} For $w_1, w_2 \in \Bbb N ,$ we have
$$\frac{1}{w_1}\sum_{i=0}^{w_1-1}\int_{\Bbb
Z_p}e^{(x_1+w_2x+\frac{w_2}{w_1}i)tw_1}dx_1
=\frac{1}{w_2}\sum_{i=0}^{w_2-1}\int_{\Bbb
Z_p}e^{(x_2+w_1x+\frac{w_1}{w_2}i)tw_2}dx_2 .$$
\endproclaim

By (18) and (18-1), we obtain the following corollary:

\proclaim{Corollary 4} For $n\geq 0$, we have
$$\sum_{i=0}^{w_1-1}B_n(w_2x+\frac{w_2}{w_1}i)w_1^{n-1}
=\sum_{i=0}^{w_2-1}B_n(w_1x+\frac{w_1}{w_2}i)w_2^{n-1}.$$
\endproclaim
Remark. Setting $w_2=1$ in Corollary 4, we get the multiplication
theorem for the Bernoulli polynomials as follows (see [26, 27]):
$$B_n(w_1x)=w_1^{n-1}\sum_{i=0}^{w_1-1}B_n(x+\frac{i}{w_1}).$$

\vskip 20pt

{\bf\centerline {\S 3. Symmetry  fermionic $p$-adic invariant
integral on $\Bbb Z_p$}} \vskip 10pt

In this section we assume that $p$ is  a fixed odd prime number. For
$f\in UD(\Bbb Z_p),$ the fermionic $p$-adic $q$-integral is defined
as
$$I_{-q}(f)=\int_{\Bbb Z_p}f(x)d\mu_{-q}(x)=\lim_{N\rightarrow
\infty}\frac{1+q}{1+q^{p^N}}\sum_{x=0}^{p^N-1}f(x)(-q)^x, \text{
see([10 ] ). }\tag19$$ Now we consider the fermionic $p$-adic
invariant integral on $\Bbb Z_p$ as
$$I_{-1}(f)= \lim_{q\rightarrow 1}I_{-q}(f)=\int_{\Bbb Z_p}f(x)d\mu_{-1}(x).
\tag20$$ From (20) we can derive the equation of fermionic $p$-adic
invariant integral on $\Bbb Z_p$ as follows:
$$\int_{\Bbb Z_p}f(x+1)d\mu_{-1}(x)=-\int_{\Bbb
Z_p}f(x)d\mu_{-1}(x)+2f(0). \tag21$$ From (21), we also note that
$$\int_{\Bbb Z_p}f(x+n)d\mu_{-1}(x)+(-1)^{n-1}\int_{\Bbb
Z_p}f(x)d\mu_{-1}(x)=2\sum_{l=0}^{n-1}(-1)^{n-1-l}f(l), \text{ for
$n\in\Bbb N $ }. \tag22$$ If we take $n(=odd)\in\Bbb N$ in (22),
then we have
$$\int_{\Bbb Z_p}f(x+n)d\mu_{-1}(x)+\int_{\Bbb
Z_p}f(x)d\mu_{-1}(x)=2\sum_{l=0}^{n-1}(-1)^lf(l). \tag23$$ By using
(23), we obtain the following equation:
$$\int_{\Bbb Z_p}e^{(x+n)t}d\mu_{-1}(x)+\int_{\Bbb
Z_{p}}e^{xt}d\mu_{-1}(x) =2\sum_{l=0}^{n-1}(-1)^le^{lt}. \tag24$$
Let $T_k(n)=\sum_{l=0}^n(-1)^ll^k$. Then $T_k(n)$ is called by the
the alternating sums of powers consecutive integers. From the
definition of the fermionic $p$-adic invariant integral on $\Bbb
Z_p$, it is not difficult to show that
$$\int_{\Bbb Z_p}e^{(x+n)t}d\mu_{-1}(x)+\int_{\Bbb
Z_p}e^{xt}d\mu_{-1}(x)=\frac{2\int_{\Bbb
Z_p}e^{xt}d\mu_{-1}(x)}{\int_{\Bbb Z_p}e^{nxt}d\mu_{-1}(x)}.
\tag25$$ In (25), we note that
$$\int_{\Bbb Z_p}e^{nxt}d\mu_{-1}(x)=\frac{2}{e^{nt}+1}.$$
Let $w_1, w_2 (\in\Bbb N)$ be odd. By using double fermionic
$p$-adic invariant integral on $\Bbb Z_p$, we have the the
following:
$$\frac{\int_{\Bbb Z_p}\int_{\Bbb
Z_p}e^{(w_1x_1+w_2x_2)t}d\mu_{-1}(x_1)d\mu_{-1}(x_2)}{\int_{\Bbb
Z_p}e^{w_1w_2xt}d\mu_{-1}(x)}=\frac{2(e^{w_1w_2t}+1)}{(e^{w_1t}+1)(e^{w_2t}+1)}.
$$
Now, we also consider the following fermionic $p$-adic invariant
integral on $\Bbb Z_p$ associated with Euler polynomials.
$$I^*= \frac{\int_{\Bbb Z_p}\int_{\Bbb
Z_p}e^{(w_1x_1+w_2x_2+w_1w_2x)t}d\mu_{-1}(x_1)d\mu_{-1}(x_2)}{\int_{\Bbb
Z_p}e^{w_1w_2xt}d\mu_{-1}(x)}=\frac{2e^{w_1w_2xt}(e^{w_1w_2t}+1)}{(e^{w_1t}+1)(e^{w_2t}+1)}.\tag26$$
From (24) and (25), we note that
$$ \frac{2\int_{\Bbb
Z_p}e^{xt}d\mu_{-1}(x)}{\int_{\Bbb
Z_p}e^{w_1xt}d\mu_{-1}(x)}=2\sum_{l=0}^{w_1-1}(-1)^le^{lt}
=2\sum_{k=0}^{\infty}\sum_{l=0}^{w_1-1}(-1)^ll^k\frac{t^k}{k!}=2\sum_{k=0}^{\infty}T_k(w_1-1)\frac{t^k}{k!}.\tag27$$
By (26) and (27) we see that
$$\aligned
I^*&=\left(\frac{1}{2}\int_{\Bbb
Z_p}e^{w_1(x_1+w_2x)t}d\mu_{-1}(x_1)\right)\left(\frac{2\int_{\Bbb
Z_p}e^{w_2x_2t}d\mu_{-1}(x_2)}{\int_{\Bbb
Z_p}e^{w_1w_2xt}d\mu_{-1}(x)}\right)\\
&=\frac{1}{2}\left(\sum_{i=0}^{\infty}E_{i}(w_{2}x)\frac{w_1^i}{i!}t^i\right)
\left(2\sum_{l=0}^{\infty}T_l(w_1-1)\frac{w_2^lt^l}{l!}\right)\\
&=\sum_{n=0}^{\infty}\left(\sum_{i=0}^n\binom{n}{i}E_i(w_2x)T_{n-i}(w_1-1)w_1^iw_2^{n-i}\right)
\frac{t^n}{n!},
\endaligned\tag28$$
where $E_n(x)$ are the $n$-th ordinary Euler polynomials.

On the other hand,
$$\aligned
I^*&=\left(\frac{1}{2}\int_{\Bbb
Z_p}e^{w_2(x_2+w_1x)t}d\mu_{-1}(x_2)\right)\left(\frac{2\int_{\Bbb
Z_p}e^{w_1x_1t}d\mu_{-1}(x_1)}{\int_{\Bbb
Z_p}e^{w_1w_2xt}d\mu_{-1}(x)}\right)\\
&=\frac{1}{2}\left(\sum_{i=0}^{\infty}E_{i}(w_{1}x)\frac{w_2^i}{i!}t^i\right)
\left(2\sum_{l=0}^{\infty}T_l(w_2-1)\frac{w_1^lt^l}{l!}\right)\\
&=\sum_{n=0}^{\infty}\left(\sum_{i=0}^n\binom{n}{i}E_i(w_1x)T_{n-i}(w_2-1)w_2^iw_1^{n-i}\right)
\frac{t^n}{n!}.\endaligned\tag29$$

By comparing the coefficients on both sides of (28) and (29),  we
obtain the following theorem.

\proclaim{Theorem 5} Let $w_1, w_2(\in \Bbb N)$ be odd and let
$n\geq 0$. Then we have
$$ \sum_{i=0}^n\binom{n}{i}E_i(w_2x)T_{n-i}(w_1-1)w_1^iw_2^{n-i}
=\sum_{i=0}^n\binom{n}{i}E_i(w_1x)T_{n-i}(w_2-1)w_2^iw_1^{n-i},$$
where $E_n(x)$ are the $n$-th ordinary Euler polynomials.
\endproclaim
Remark. Setting $x=0$ in Theorem 5 we obtain
$$\sum_{i=0}^n\binom{n}{i}E_iT_{n-i}(w_1-1)w_1^iw_2^{n-i}
=\sum_{i=0}^n\binom{n}{i}E_iT_{n-i}(w_2-1)w_2^iw_1^{n-i},$$ where
$E_n$ are the $n$-th ordinary Euler numbers. If we take $w_2=1$ in
Theorem 5, then we have
$$E_n(w_1x)=\sum_{i=0}^n\binom{n}{i}E_i(x)T_{n-i}(w_1-1)w_1^i.
\tag30$$ Setting $x=0$ in (30) we obtain the following corollary.

\proclaim{ Corollary 6} Let $w_1(\in \Bbb N)$ be odd and let
$n\in\Bbb N$. Then we have
$$E_n=\left(\frac{1}{1-w_1^n}\right)\sum_{i=0}^{n-1}\binom{n}{i}E_iT_{n-i}(w_1-1)w_1^i.$$
\endproclaim
From (27), we note that
$$\aligned
I^*&=\left(\frac{e^{w_1w_2xt}}{2}\int_{\Bbb
Z_p}e^{w_1x_1t}d\mu_{-1}(x_1)\right)\left(\frac{2\int_{\Bbb
Z_p}e^{w_2x_2t}d\mu_{-1}(x_2)}{\int_{\Bbb
Z_{p}}e^{w_1w_2xt}d\mu_{-1}(x)} \right)\\
&=\left(\frac{e^{w_1w_2xt}}{2}\int_{\Bbb
Z_p}e^{w_1x_1t}d\mu_{-1}(x_1)\right)\left(2\sum_{l=0}^{w_1-1}(-1)^le^{w_2lt}\right)\\
&=\sum_{l=0}^{w_1-1}(-1)^l\int_{\Bbb
Z_p}e^{(x_1+w_2x+\frac{w_2}{w_1}l)tw_1}d\mu_{-1}(x_1)\\
&=\sum_{n=0}^{\infty}\left(w_1^n\sum_{l=0}^{w_1-1}(-1)^lE_n(w_2x+\frac{w_2}{w_1}l)\right)
\frac{t^n}{n!}.
\endaligned\tag31$$
On the other hand,
$$\aligned
I^*&=\left(\frac{e^{w_1w_2xt}}{2}\int_{\Bbb
Z_p}e^{w_2x_2t}d\mu_{-1}(x_2)\right)\left(\frac{2\int_{\Bbb
Z_p}e^{w_1x_1t}d\mu_{-1}(x_1)}{\int_{\Bbb
Z_{p}}e^{w_1w_2xt}d\mu_{-1}(x)} \right)\\
&=\left(\frac{1}{2}\int_{\Bbb
Z_p}e^{w_2x_2t}d\mu_{-1}(x_2)\right)\left(2\sum_{l=0}^{w_{2}-1}(-1)^le^{(w_1x +\frac{w_1}{w_2}l)w_2t}\right)\\
&=\sum_{l=0}^{w_2-1}(-1)^l\int_{\Bbb
Z_p}e^{(x_2+w_1x+\frac{w_1}{w_2}l)tw_2}d\mu_{-1}(x_2)\\
&=\sum_{n=0}^{\infty}\left(w_2^n\sum_{l=0}^{w_2-1}(-1)^lE_n(w_1x+\frac{w_1}{w_2}l)\right)
\frac{t^n}{n!}.
\endaligned\tag32$$
By comparing the coefficients on the both sides of (31) and (32), we
obtain the following theorem.

\proclaim{ Theorem 7} Let $w_1, w_2 (\in \Bbb N)$ be odd and let
$n\geq 0$. Then we have
$$w_1^n\sum_{l=0}^{w_1-1}(-1)^lE_n(w_2x+\frac{w_2}{w_1}l)
=w_2^n\sum_{l=0}^{w_2-1}(-1)^lE_n(w_1x+\frac{w_1}{w_2}l).$$
\endproclaim

Setting $w_2=1$ in Theorem 7, we get the multiplication theorem for
the Euler polynomials as follows:
$$E_n(w_1x)=w_1^n\sum_{l=0}^{w_1-1}(-1)^lE_n(x+\frac{l}{w_1}).$$

\vskip 20pt

\quad Taekyun Kim

\quad EECS, Kyungpook National University, Taegu 702-701, S. Korea

\quad e-mail:\text{ tkim$\@$knu.ac.kr; tkim64$\@$hanmail.net}

\enddocument